\documentclass{article}

\usepackage{arxiv}

\usepackage{bm}
\usepackage{subcaption}
\usepackage{amsmath}
\usepackage{amsfonts}
\usepackage{listings}

\usepackage[svgnames]{xcolor}

\lstdefinelanguage{XML}
{
  basicstyle=\ttfamily\footnotesize,
  morestring=[b]",
  moredelim=[s][\bfseries\color{Maroon}]{<}{\ },
  moredelim=[s][\bfseries\color{Maroon}]{</}{>},
  moredelim=[l][\bfseries\color{Maroon}]{/>},
  moredelim=[l][\bfseries\color{Maroon}]{>},
  morecomment=[s]{<?}{?>},
  morecomment=[s]{<!--}{-->},
  commentstyle=\color{DarkOliveGreen},
  stringstyle=\color{blue},
  identifierstyle=\color{red}
}

\newcommand{\ra}{{\rm a}}

\usepackage{color}

\usepackage{pgfplots}
\usepackage{pgfplotstable}
\pgfplotsset{
compat=newest, 
tick label style={font=\footnotesize}, 
}

\usepackage{overpic}

\usepackage{tikz}
\usetikzlibrary{mindmap,shadows}
\usetikzlibrary{arrows,intersections}
\usepackage{smartdiagram}

\usepackage{booktabs}
\usepackage[]{units}

\newtheorem{theorem}{Equation}[section]

\newcommand{\rd}{{\,\rm d}}

\title{openCFS: Open Source Finite Element Software for Coupled Field Simulation - Part Acoustics}

\author{ Stefan Schoder \\
	Group of Aeroacoustics and Vibroacoustics, IGTE\\
	TU Graz\\
	Inffeldgasse 16c, 8010 Graz \\
	\texttt{stefan.schoder@tugraz.at} \\
	\And
	Klaus Roppert \\
	Group of Multiphysics, IGTE\\
	TU Graz\\
	Inffeldgasse 18, 8010 Graz \\
	\texttt{klaus.roppert@tugraz.at} \\
}

\renewcommand{\shorttitle}{openCFS Software}

\usepackage{hyperref}
\hypersetup{
pdftitle={openCFS Software},
pdfsubject={preprint},
pdfauthor={Schoder, Roppert},
pdfkeywords={FEM Sotware, open Source},
}

\begin{document}
\maketitle

\begin{abstract}
	The finite element method offers attractive methods for the numerical solution of coupled field problems arising in sensors and actuator simulations of various physical domains, like electrodynamics, mechanics, and thermodynamics. 
    With this application perspective and being open, accessible, and fast implementations are possible, openCFS was launched in 2020. It provides an open-source framework for implementing partial differential equations using the finite element method. In particular, the acoustic module is part of active development, including several key methods. These methods include the perfectly-matched layer technique, non-confirming interface formulations, Lagrangian basis function, Legendre basis functions, spectral element formulations, a nodal element type, edge-based element type (aeroacoustic post-processing), absorbing boundary conditions, frequency dependent-material for time-harmonic and time-dependent simulations. Time-dependent simulations, time-harmonic simulations, and eigenvalue simulations are supported. Several variants of acoustic equations are implemented, including the relevant source terms and wave operators for aeroacoustics. The package includes rotating domains and non-conforming interfaces for fan noise simulations. It also contains an API to the Python3 package pyCFS \cite{wurzinger2024pycfs}. This paper presents openCFS with a focus on the acoustic module.
\end{abstract}

\keywords{Open Source FEM Software \and Multiphysics Simulation \and C++ \and Acoustics \and Aero-Acoustics \and openCFS}

%

\section{Introduction}
\label{sec:Intro} 

The finite element method (FEM) progressed significantly in solving wave-based acoustics around 2000 \cite{harari2006survey,zienkiewicz2000achievements,colton2000recent,ihlenburg1998finite}. The origins of openCFS started at that time. In 2020, CFS++ went open source (openCFS) and provides an object-oriented finite element framework for research in coupled-field simulations using FEM. Besides these standard FEM, many different numerical techniques to solve partial differential equations were established and found to provide reliable solutions \cite{dumbser2007arbitrary,fiala2014nihu,tam1993dispersion}. The overarching aim of the methods is to derive a numerical scheme suitable for yielding the solution of a partial differential equation. Regarding the different numerical methods, the standard Bubnov–Galerkin method was chosen to establish the weak formulation for the respective coupled field implementations. 

Its flexible meshing and interfacing between different field equations using coupling terms is the core strength of openCFS. Over the last years, the module acoustic (\textit{openCFS-Acoustics}) and the coupling to the mechanic field (vibroacoustic \cite{engelmann2020generic}, with examples see\footnote{\url{https://gitlab.com/openCFS/Testsuite/-/tree/master/TESTSUIT/Coupledfield/MechAcou}}) and fluid dynamics as a source term (aeroacoustics \cite{schoder2024aeroacoustic,schoder2023acoustic,rucz2024analysis}) was one of the foci. For instance, vibroacoustic simulations are performed using full coupling inside openCFS, solving the acoustics and solid mechanical field in parallel \cite{tieber2024hybrid}, and the coupling is realized via coupling matrices (this idea is described in \cite{kaltenbacher2015numerical}). Another possibility is to prescribe the surface normal velocity as a boundary condition and perform a forward coupling simulation using the acoustic module only \cite{wurzinger2024experimental}; this is only recommended in special occasions where feedback is negligible (e.g., sound radiation from a stiff, closed mechanical structure). In contrast, some cases exist where the acoustics have strong back-coupling properties onto the solid mechanical field \cite{nager2023investigation,kraxberger2023alignment}. In aeroacoustics, usually, the segregated workflow is applied \cite{schoder2023offline} where in a separate flow simulation using a state-of-the-art flow solver is calculating the flow solution \cite{schoder2019hybrid} and then integrated by the module openCFS-Data \cite{schoder2023opencfs} onto the acoustic grid as aeroacoustic source terms that are then propagated using an acoustic simulation. The particular focus of interest is numerical methods to solve these coupled fields and connected applications. The convergence of simulations using non-conforming interfaces for stationary and rotating domains \cite{Hueppe,kaltenbacher2018non,schoder2018aeroacoustic,roppert2018non} (e.g. example with rotating interface\footnote{\url{https://gitlab.com/openCFS/Testsuite/-/tree/master/TESTSUIT/Singlefield/Acoustics/RotationNCAcouPotential}}) have been performed in \cite{schoder2021application}. The particular strength is the strong connection to engineering research and the regular use of the software for source code and software solution development in industrial and basic research projects. It is hosted on GitLab (gitlab.com/openCFS), with the accompanying package API pyCFS, the automated build environment, the testing examples, and the documentation. The software includes several physical fields
\begin{itemize}
\item Acoustics (\textit{openCFS-Acoustics})
\item Electrodynamics (\textit{openCFS-Edyn})
\item Mechanics (\textit{openCFS-Mechanics})
\item Piezoelectrics (\textit{openCFS-Piezo})
\item Heat transfer (\textit{openCFS-Heat}).
\end{itemize}

The paper is organized as follows: In Sec.~\ref{Sec:FEM}, we provide the generic finite element formulation of the standard wave equation as implemented in openCFS.
In Sec. \ref{Sec:W}, \ref{Sec:FEM}, and \ref{Sec:A}. We outline how to implement a finite element formulation in the framework of openCFS. We start by stating the wave equation and deriving its weak formulation for finite element implementation. Then, we discuss the calling structure of openCFS and show the implementation. We show the basic building blocks of a finite element simulation using openCFS. The motivation of this publication is to encourage users to use and develop openCFS with us (the usage is free under the MIT license). Section \ref{Sec:Applications} discusses selected publications and concludes the article.

\section{Wave equation} \label{Sec:W}

\subsection{Standard wave equation}
\begin{figure}[htp!]
\centering
\includegraphics[scale=0.8]{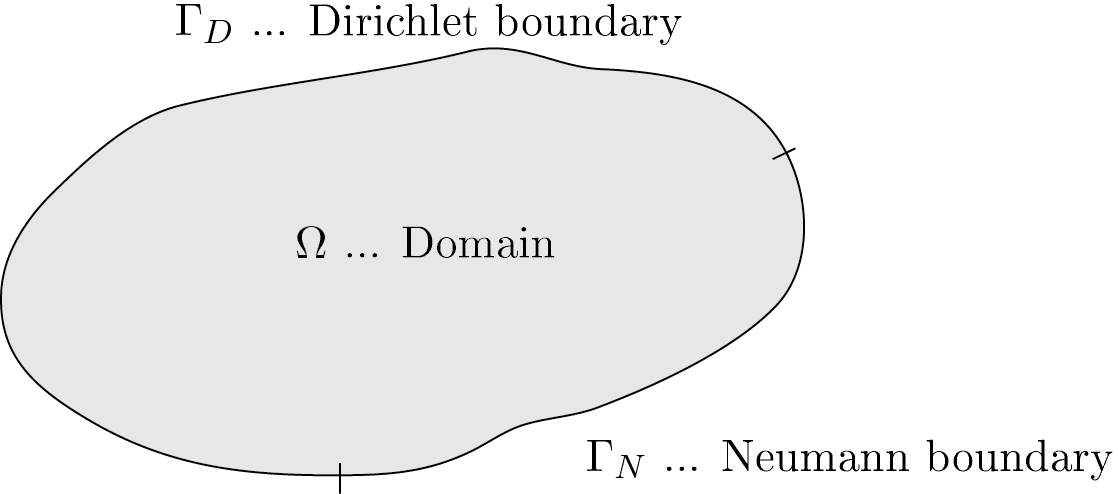}
\caption{\label{fig:domain} Domain of the wave equation.}
\end{figure}

A basic description of the FEM for acoustics in connection with openCFS can be found in \cite{kaltenbacher2015numerical}. 
\begin{theorem}[Wave equation]
Given $c$, $p_0^\ra(\bm x)$, $\rd p_0^\ra(\bm x)$, $p_{\Gamma_D}^\ra(\bm x,t)|_{\Gamma_D}$, $u_{n\Gamma_N}^\ra(\bm x,t)|_{\Gamma_N}$, domain $\bm x \in \Omega$, and time $t \in [0,T]$ find $p_\ra$ such that
\begin{enumerate}
\item Partial differential equation:
\[
\frac{1}{c^2} \frac{\partial^2 p_\ra}{\partial t^2} 
 -  \nabla \cdot \nabla p_\ra = f \,\,\,\,\,\, \text{in} \,\,  \Omega \times [0,T]
\]
\item Boundary conditions:
\[
 p_\ra(\bm x,t) = p_{\Gamma_D}^\ra(\bm x,t) \,\,\,\,\,\, \text{at} \,\,  \Gamma_D \times [0,T]
\]
\[
\nabla p_\ra(\bm x,t)\cdot \bm n = u_{n\Gamma_N}^\ra(\bm x,t) \,\,\,\,\,\, \text{at} \,\,  \Gamma_N \times [0,T]
\]
\item Initial conditions:
\[
 p_\ra(\bm x,0) =p_0^\ra(\bm x) \,\,\,\,\,\, \forall \bm x \in \Omega
\]
\[
\frac{\partial}{\partial t} p_\ra(\bm x,0) = \dot{p}_\ra(\bm x,0) =\rd p_0^\ra(\bm x) \,\,\,\,\,\, \forall \bm x \in \Omega
\]
\end{enumerate}
\end{theorem}
Regarding the definition of the wave equation on the domain $\Omega$, the Dirichlet boundary $\Gamma_\mathrm{D}$, and the Neumann boundary $\Gamma_\mathrm{N}$ (see Fig. \ref{fig:domain}), we derive the weak formulation of the wave equation. Examples can be found here\footnote{\url{https://gitlab.com/openCFS/Testsuite/-/tree/master/TESTSUIT/Singlefield/Acoustics}} and for a first example we recommend\footnote{\url{https://gitlab.com/openCFS/Testsuite/-/tree/master/TESTSUIT/Singlefield/Acoustics/Abc2d}}.

\subsection{Alternative wave equations implemented}

Besides the standard wave equation, several variants of wave equations are implemented in openCFS.

\subsubsection{Frequency dependent speed of sound (time-harmonically)}
\newcommand{\eqVar}[1] {#1_{\mathrm{eq}}}
\newcommand{\acVar}[1] {#1_{\mathrm{a}}}
The modified Helmholtz equation with frequency-dependent material reads as
\begin{equation}
\label{eq:inhomogWaveEq}
\frac{\omega^2}{\eqVar{K}(\omega)}  \acVar{p}
+ \nabla \cdot \left( \frac{1}{\eqVar{\rho}(\omega)}  \nabla \acVar{p} \right)  = 0 \, ,
\end{equation}
with a frequency-dependent equivalent bulk modulus $\eqVar{K}(\omega)$ and density $\eqVar{\rho}(\omega)$. This type of wave equation is used to model various absorber materials for acoustic treatment \cite{czwielong2023mikroperforierte,floss2019finite}. An example simulation for an acoustic meta-material placed in a duct can be found in the openCFS-Testsuite \footnote{\url{https://gitlab.com/openCFS/Testsuite/-/tree/master/TESTSUIT/Singlefield/Acoustics/Muffler1MPPnmg}}.

\subsubsection{Frequency dependent speed of sound (time-domain)}
Using the inverse Fourier transform eq. \eqref{eq:inhomogWaveEq}  yields in time domain
\begin{equation}
\{\frac{1}{\eqVar{K}} * \acVar{\ddot{p}}\}(t)
- \{\nabla \cdot \left( \frac{1}{\eqVar{\rho}} * \nabla \acVar{p} \right)\}(t)  = 0 \, ,
\end{equation}
with $\acVar{p}$ denoting the acoustic pressure, $\{a*b\}(t) = \int_{-\infty}^{\infty} a(t-\tau)\cdot b(\tau)\ d\tau$ the time-convolution operator, and the two dots above the variable indicating the second-order partial time derivative. This equation can be used only when time-dependent simulations, such as turbomachinery simulations, are necessary. An example simulation for an acoustic meta-material can be found in the openCFS-Testsuite \footnote{\url{https://gitlab.com/openCFS/Testsuite/-/tree/acou_TimeDomainEquivalentFluid_New/TESTSUIT/Singlefield/Acoustics/TDEFMultiMaterial2D}} and a non-reflecting boundary condition here\footnote{\url{https://gitlab.com/openCFS/Testsuite/-/tree/acou_TimeDomainEquivalentFluid_New/TESTSUIT/Singlefield/Acoustics/TDEFwNRBC2D}}.

\subsubsection{Spatially dependent speed of sound}

A space-dependent speed of sounds, varying for every element in the domain, can be prescribed. A temperature field governs the dependency and could model the refraction effects of the speed of sound variability in atmospheres or a heated jet. 

\begin{equation}
  \acVar{\ddot{p}}
- \nabla \cdot \left( c^2(T)  \nabla \acVar{p} \right)  = 0 \, .
\end{equation}
An example simulation with a temperature-dependent PML can be found in the openCFS-Testsuite \footnote{\url{https://gitlab.com/openCFS/Testsuite/-/blob/master/TESTSUIT/Singlefield/Acoustics/WaveFlowTempPml}}.

\subsubsection{Convective wave equation}
Presented by Pierce \cite{pierce1990wave}, a convective wave equation in the acoustic velocity potential $\acVar{\bm u} = -\nabla \psi$ with a spatially uniform speed of sound and a mean flow velocity $\overline{\boldsymbol{u}}$ is given by 
$$\frac{1}{c_0^2} \frac{\mathrm{D}^2 \psi}{\mathrm{D} t^2}-\nabla \cdot \nabla \psi=0$$
with the substantial time derivative
$$\frac{\mathrm{D}}{\mathrm{D} t}=\frac{\partial}{\partial t}+\overline{\boldsymbol{u}} \cdot \nabla\, .$$ 
An example simulation and usage can be found in the openCFS-Testsuite \footnote{\url{https://gitlab.com/openCFS/Testsuite/-/tree/master/TESTSUIT/Singlefield/Acoustics/PCWE}}.

\subsection{Pressure/Velocity form of the wave equation}
The advection form of the wave equation in terms of the pressure and particle velocity, as unknown, is provided in the AcousticMixedPDE method. Today, no active development is carried out here. Example simulations and usage can be found in the openCFS-Testsuite \footnote{\url{https://gitlab.com/openCFS/Testsuite/-/tree/master/TESTSUIT/Singlefield/AcousticMixed}}.

\subsubsection{Viscous acoustics (Linearized Navier-Stokes Equation)}
In the method LinFlowPDE, the linearized Navier-Stokes equations are implemented. For example, one can look at \cite{guilvaiee2023validated,mayrhofer2024moving}. Example simulations and usage can be found in the openCFS-Testsuite \footnote{\url{https://gitlab.com/openCFS/Testsuite/-/tree/master/TESTSUIT/Singlefield/LinFlow}}.

\subsection{Postprocessing of acoustic fields (acousticSplitPDE)}
A general compressible fluid field can be decomposed into a longitudinal and transversal process using the acousticSplitPDE \cite{schoder2020postprocessing,schoder2020postprocessing2,schoder2019helmholtz,schoder2023acoustic,schoder2019revisiting}. Example simulations and usage can be found in the openCFS-Testsuite \footnote{\url{https://gitlab.com/openCFS/Testsuite/-/tree/master/TESTSUIT/Singlefield/Split}}.

\section{Finite Element Formulation} \label{Sec:FEM}

\subsection{Standard wave equation}
We define the solution function space for the acoustic pressure $
 p_\ra \in \mathcal{V}=\lbrace v \in H^1(\Omega) | v =p_0^\ra \,\,  \text{at} \,\,  \Gamma_D  \rbrace
$
and the test function space $
w \in \mathcal{W}=\lbrace u \in H^1(\Omega) | u =0 \,\,  \text{at} \,\,  \Gamma_D  \rbrace
$. Using this function space, we can obtain the weak formulation of the wave equation. Additionally, we restrict the function spaces to a discrete (limited number of functions) solution function space $
 p_\ra^h \in \mathcal{V}^h \subset \mathcal{V}; \,\, \mathcal{V}^h = \lbrace v^h \in H^1(\Omega) | v^h =p_0^\ra \,\,  \text{at} \,\,  \Gamma_D  \rbrace
$ and a discrete test function space $
 w^h \in \mathcal{W}^h \subset \mathcal{W}; \,\, \mathcal{W}^h = \lbrace u^h \in H^1(\Omega) | u^h =0 \,\,  \text{at} \,\,  \Gamma_D  \rbrace
$. In doing so, we obtain the semi-discrete finite element form.
\begin{theorem}[Weak formulation (semi-discrete)]
Given $c$, $p_0^\ra(\bm x)$, $\rd p_0^\ra(\bm x)$, $p_{\Gamma_D}^\ra(\bm x,t)|_{\Gamma_D}$, $u_{n\Gamma_N}^\ra(\bm x,t)|_{\Gamma_N}$, domain $\bm x \in \Omega$, and time $t \in [0,T]$. Find $p_\ra^h \in \mathcal{V}^h$ such that $\forall$ $w^h \in \mathcal{W}^h$ the equation holds:
\begin{enumerate}
\item Weak form:
\begin{eqnarray*}
\int\limits_\Omega \frac{1}{c^2}\, w^h\, \ddot p_\ra^h \, d\Omega
+ \int\limits_\Omega\, \nabla w^h\, \cdot \, \nabla p_\ra^h\, 
d\Omega -  \int\limits_{\Gamma_N}\,w^h u_{n\Gamma_N}^\ra\, 
ds
 &=& \int\limits_\Omega \, w^h f \, d\Omega
\end{eqnarray*}
\item Boundary conditions:
\[
 p_\ra(\bm x,t) = p_{\Gamma_D}^\ra(\bm x,t) \,\,\,\,\,\, \text{at} \,\,  \Gamma_D \times [0,T]
\]
\item Initial conditions:
\[
 p_\ra(\bm x,0) =p_0^\ra(\bm x) \,\,\,\,\,\, \forall \bm x \in \Omega
\]
\[
\frac{\partial}{\partial t} p_\ra(\bm x,0) = \dot{p}_\ra(\bm x,0) =\rd p_0^\ra(\bm x) \,\,\,\,\,\, \forall \bm x \in \Omega
\]
\end{enumerate}
\end{theorem}
Inserting the finite element ansatz ($n_{\rm n}$ number of global basis functions)
\begin{eqnarray*}
p_\ra &\approx & p_\ra^h =  \sum\limits_{b=1}^{n_{\rm n}} \, 
N_b(\bm x) p_{\ra,b}(t) \,; \\
w &\approx&  w^h =  \sum\limits_{c=1}^{n_{\rm n}} \, N_c(\bm x) w_c(t) \,;
\end{eqnarray*}
The basis functions are known. The coefficients are the unknowns and model the physical solution. For the FE ansatz, higher order finite element basis functions \cite{hueppe2012spectral} can be used. The semi-discrete form (space discrete, time continuous) can be expressed by an algebraic system of equations
\[
{\bf M} \underline{\ddot{p}_\ra}^h +
{\bf K} \underline{p_\ra}^h = \underline f^h + \underline f^h_{Neumann} + \underline f^h_{Dirichlet} \, .
\]
The underline acoustic pressure $\underline{p_\ra}^h$ denotes an algebraic vector collecting the unknowns of the scalar acoustic pressure. $\underline{f}^h$, $\underline f^h_{Neumann}$, $\underline f^h_{Dirichlet}$ denote the right-hand sides according to a given source term, the inhomogeneous Neumann boundary conditions, and the inhomogeneous Dirichlet boundary condition, respectively. The mass and stiffness matrice are given by $\boldsymbol{M}, \boldsymbol{K} \in \mathbb{R}^{n_{\mathrm{eq }}} \times \mathbb{R}^{n_{\mathrm{eq}}}$ with $n_{\mathrm{eq }}$ the number of unknowns. $n_{\mathrm{eq }} < n_{\rm n}$ and the discrepancy is given by the elimination routine when integrating the boundary conditions. An illustrative explanation of this can be given by the null space of the system \cite{schoder2024approach}. The finite difference time discretization can be performed by the implicit Hilber-Hughes-Taylor (HHT) method \cite{hilber1977improved}. It recovers a classical Newmark method under certain parameters. It is second-order temporal accurate and allows for energy dissipation, and under certain parameters, it is unconditionally stable.

In openCFS, the implementation follows a basic work principle. Firstly, the constructor needs to specify some basic variables and flags, like the material class type, nonlinearity, couplings, geometry updates during the simulation, and, in addition, some specifics of the acoustic FEM implementation. openCFS uses a physical field-driven approach and terminology, and the respective formulations are collected by their physical domain (like acoustics). This means that different weak formulations that belong to the scientific field of acoustics are collected, with all benefits and consequences. In the acousticPDE (this is the name of the main file of the acoustic FEM implementation in openCFS), the following flavors of acoustic variables (acoustic pressure, acoustic potential) can be found. The possible analysis types (e.g., time-harmonic) are also specified here. Subsequently, the operators (integrators) are defined as bilinear forms.
These bilinear forms form the system matrices $\bf K, \bf M$ and are connected to the discretization of the continuous bilinear form, the basis functions, the element definitions, the FE mesh, the continuous and discrete operator, and the assembly routine. Dirichlet boundaries are handled by elimination. The speed of sound is defined by a constant or frequency-dependent density and bulk modulus \cite{kraxberger2023validated}. Modeling a frequency-dependent material in the time domain can be handled by an auxiliary variable technique \cite{maurerlehner2023time,schoder2024metamat}; a similar approach was conducted for the perfectly matched layer technique in the time domain \cite{kaltenbacher2013modified}. A linear form implements the right-hand side, and a continuous source density function can be applied as forcing. Alternatively, such sources can be read from an external file (e.g., Synthetic noise generation and radiation \cite{schoder2023opencfsSNGR}). Rayleigh damping can be used in the acoustic domain as well. 
In addition to the Neumann and Dirichlet boundary conditions, an impedance and boundary layer condition can be described at the boundaries. Furthermore, a solve step is defined and executed by the time integration routine (finite difference method). At the end of each timestep, variables defined in the post-processing results are executed. The acousticPDE is a subclass of the singlePDE, where the singlePDE describes the basic class of a single field FEM routine and is a subclass of the stdPDE, which is a subclass of the basePDE. Each PDE method must at least implement the methods CreateFeSpaces (finite element space of the unknown physical quantity), DefineIntegrators (the bilinear forms of the weak form), DefineNcIntegrators (includes the handling of the non-conforming mesh interfaces), DefineSurfaceIntegrators (defines bilinear forms at the surface, for instance, the absorbing boundary condition), DefineRhsLoadIntegrators (excitation of the PDE), InitTimeStepping (definition of the time stepping scheme), DefineSolveStep (defines the step in which the equation is solved), DefinePrimaryResults (defines the primary quantity as a result of the simulation), DefinePostProcResults (defines additional results of a simulation, for instance, the acoustic intensity).


\subsection{Convective wave equation}

Similarly, the weak formulation of the convective wave equation can be obtained \cite{kaltenbacher2021stable}
$$
\begin{array}{r}
\frac{1}{c_0^2} \int_{\Omega} w \ddot{\psi} \mathrm{~d} \boldsymbol{x}+\frac{1}{c_0^2} \int_{\Omega} w(\overline{\boldsymbol{u}} \cdot \nabla) \dot{\psi} \mathrm{~d} \boldsymbol{x}-\frac{1}{c_0^2} \int_{\Omega} \dot{\psi}(\overline{\boldsymbol{u}} \cdot \nabla) w \mathrm{d} \boldsymbol{x} \\
-\frac{1}{c_0^2} \int_{\Omega}(\overline{\boldsymbol{u}} \cdot \nabla) w(\overline{\boldsymbol{u}} \cdot \nabla) \psi \mathrm{d} \boldsymbol{x}-\int_{\Omega} \nabla w \cdot \nabla \psi \mathrm{d} \boldsymbol{x}+\int_{\Gamma} w \boldsymbol{n} \cdot \nabla \psi \mathrm{d} \boldsymbol{s}=0 \, . \label{eq:stableConv}
\end{array}
$$
The standard FE ansatz can be applied with appropriate FE basis functions $N_i(\boldsymbol{x})$
$$
w \approx w^h=\sum_i^{n_{\rm n}} N_i(\boldsymbol{x}) w_i(t)
$$

$$\psi \approx \psi^h=\sum_k^{n_{\rm n}} N_k(\boldsymbol{x}) \psi_k(t)\, .
$$

This will result in the semi-discrete Galerkin formulation
$$
{\bf M} \ddot{\underline{\psi}}^h+{\bf C} \dot{\underline{\psi}}^h+{\bf K} \underline{\psi}^h=\underline{f}^h .
$$
The underline acoustic potential $\underline{\psi}^h$ denotes an algebraic vector collecting the unknowns of the scalar acoustic potential. $\underline{f}^h$ denotes the right-hand side according to a given source term or boundary conditions and $\boldsymbol{C} \in \mathbb{R}^{n_{\mathrm{eq }}} \times \mathbb{R}^{n_{\mathrm{eq}}}$ is the damping matrix. The mean velocity interpolation is element-wise constant inside the convective integrators and, therefore, $L^2$. A flow list inside the simulation defines the mean velocity, which builds up the additional bilinear forms according to the stable weak form \eqref{eq:stableConv}. The major application of this equation conducted in aeroacoustics, forming the operator of the perturbed convective wave equation (PCWE) \cite{kaltenbacher2017computational,schoder2021application,tieghi2022machine,schoder2022pcwe,schoder2022cpcwe} and the aeroacoustic wave equation based on Pierce operator (AWE-PO) \cite{schoder2022aeroacoustic}. The PCWE model was used to investigate various aspects of the human phonation process \cite{schoder2020hybrid,valavsek2019application,zorner2016flow,schoder2021aeroacoustic,falk20213d,lasota2021impact,maurerlehner2021efficient,schoder2022learning}. In connection with flexible discretization based on mortar non-conforming interfaces \cite{Hueppe,kaltenbacher2018non,schoder2018aeroacoustic,roppert2018non}. These interfaces also allow for relative movement and, therefore, a handy implementation for handling rotating domains as encountered in turbomachinery simulations. Applications included axial fan noise simulations \cite{schoder2020computational,tautz2018source,kaltenbacher2017computational,kaltenbacher2016modeling,tieghi2022machine} and sound emission predictions of the turbocharger compressor \cite{kaltenbacher2020modelling,freidhager2020influences}, HVAC systems \cite{kaltenbacher2016computational,grabinger2018efficient}, and wind turbines \cite{weber2016computational}. The aeroacoustic sources are generally computed out of the flow with the auxiliary package openCFS-Data \cite{schoder2023opencfs,schoder2023opencfsSNGR,schoder2023implementation} and, for instance, presented in \cite{schoder2020aeroacoustic,schoder2020radial,schoder2019conservative}. Connected to a variety of aeroacoustic source formulations, simulations for car frame noise are noteworthy \cite{engelmann2020generic,freidhager2021simulationen,schoder2020numerical,weitz2019numerical}.

\section{Software runtime cycle of a simulation} \label{Sec:A}
The package's compilation, including external packages (solvers, modules, geometric libraries, etc.), is carried out using cmake. The compilation consists of openCFS itself, cfstool, and openCFS-Data. The openCFS C++ runtime of the finite element software executable model follows the steps of the driver class execution during a simulation. 

\subsection{Simulation setup details in XML input} \label{Sec:XML}
The simulation file of openCFS is an XML-based input format \cite{CFSdocu}. The underlying structure of our XML files is defined via a so-called XML scheme. The path to this scheme is defined in the header of every XML file, and per default, it points to our default public scheme, provided on our GitLab server:
\begin{lstlisting}[language=XML]
<?xml version="1.0"?>
<cfsSimulation ...
xsi:schemaLocation=".../CFS.xsd">
\end{lstlisting}
In general, an openCFS XML file consists of three main parts:
\begin{itemize}
\item fileFormats: definition of input-/output- and material-files
\item domain: region definitions (assign material to different regions)
\item sequenceStep: contains the definition of the analysis type and the PDE.
\end{itemize}
Several sequenceSteps can be concatenated, and results from previous sequenceSteps can be processed to produce the desired result.
\begin{lstlisting}[language=XML]
<!-- define which files are needed for simulation input & output-->
<fileFormats>
    <input>
        <!-- read mesh -->
    </input>
    <output>
        <!-- define output -->
    </output>
    <materialData file="../material/mat.xml" format="xml"/>
</fileFormats>

<domain geometryType="3d">
    <regionList>
        <!-- region definition -->
    </regionList>
</domain>

<sequenceStep index="1">
    <analysis>
        <!-- analysis type: static, transient, harmonic, multiharmonic, eigenfrequency -->
    </analysis>
    <pdeList>
        <!--for example, consider the acoustic PDE-->
        <acoustic>
            <regionList>
                <!--define on which regions the PDE is solved-->
            </regionList>
            <bcsAndLoads>
                <!--define appropriate BC's and Loads-->
            </bcsAndLoads>
            <storeResults>
                <!--define the results to be stored-->
            </storeResults>
        </acoustic>
    </pdeList>
</sequenceStep>
</cfsSimulation>
\end{lstlisting}
Detailed information about the XML-input files can be found online \cite{CFSdocu}.

\subsection{Executation cycle}
\begin{figure}[htp!]
\centering
\includegraphics[scale=0.7]{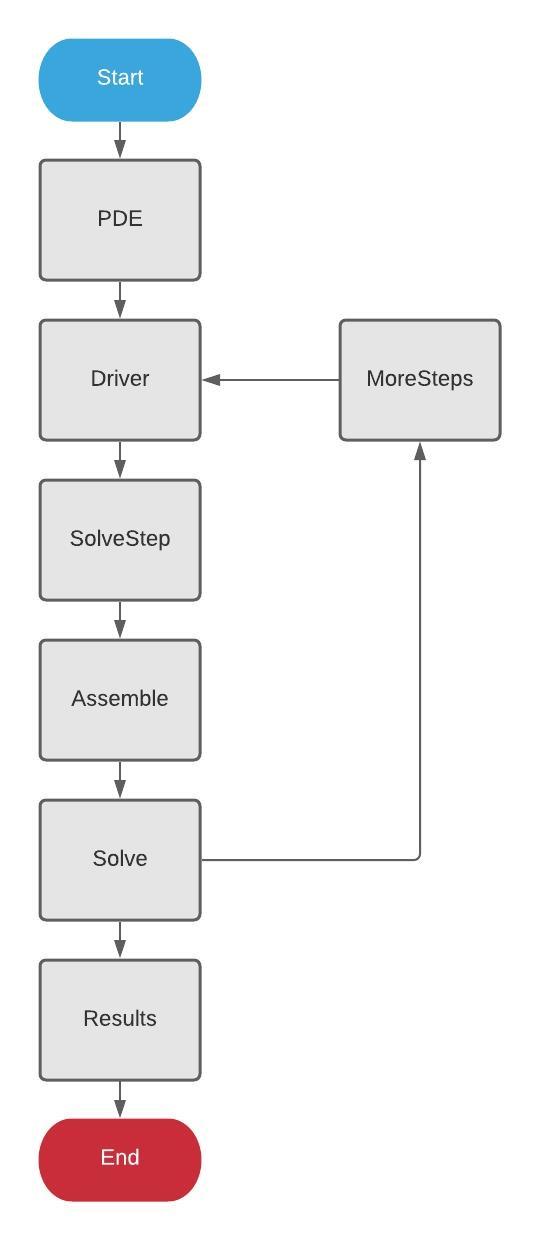}
\caption{\label{Fig:LC} openCFS life-cycle during a run.}
\end{figure}

The runtime takes the XML input file and handles the following execution steps (see Fig. \ref{Fig:LC}). It creates the partial differential equation that needs to be solved with all the setups necessary to provide a solvable algebraic system. The core simulation driver/step-iterator takes over and executes the steps of the simulation (sequence steps, time steps, frequency bins, or non-linear iterations). All in all, these can be combined in several repeated executions per simulation. As a typical finite element software, the following time discretization is possible:
\begin{itemize}
\item Time is discretized by a time-stepping scheme.
\item Time is transformed into the Fourier space (time-harmonic simulation).
\item The system of equations is transformed into an eigenvalue system.
\end{itemize} A sequence step is a step of a whole sub-simulation defined in the XML input scheme, it is intended to provide an initial field or pre-processing of an input field for a second sequence step for segregated coupled field simulations. A (solveStep) time step or frequency bin is a physical solving step encountered by the numerical procedure of time discretization. A linear static simulation has only one solve step. A non-linear iteration is part of the numerical solving procedure of the non-linear differential equation. In the solve step, the solveStep-dependent assembly routines are carried out next. 
Only the dependent parts are assembled herein; assemblies that are not dependent are (pre-)executed during the first solveStep and re-used during the additional steps. Finally, the algebraic solver takes over and solves the system of equations. Before the end of the simulation, the results are written to the files requested. Usually, when implementing a partial differential equation, no adaptions to the driver have to be made. Regarding the implementation, the following aspects have to be considered:
\begin{itemize}
\item \textbf{That is the core part you must implement. There is no way around it unless the partial differential equation is already implemented:}Bilinear forms, also called integrators, are the realizations of the left-hand side of a weak formulation and are constructed by the operators (e.g., a gradient operator). In the case of the presented wave equation, we use symmetric bilinear forms to construct the algebraic system.
    \item \textbf{Are you planning to implement something that touches the finite element space?} The finite element space manages the equation numbers for the associated unknowns. It holds and manages associated reference elements. It also knows how to incorporate boundary conditions. Furthermore, it distinguishes between geometrical elements (mesh representation) and computational elements (how the physical quantities are approximated).
\item \textbf{Are you planning to implement a new finite element basis?} The finite element function represents the function space and knows its particular selection from different possibilities.
\item \textbf{Do you consider an inhomogeneous partial differential equation?} Linear forms are integrators for the right-hand side of a partial differential equation. 
\end{itemize}

\section{Conclusion with recent application} \label{Sec:Applications}

Although many commercial simulation tools have been developed in the last years and are on the market, there is still a strong need for appropriate tools capable of simulating multi-field problems. Implementing a new partial differential equation (PDE) or physical coupling between PDEs into commercial software (via user routines) results, in most cases, in an environment that is too unflexible for forefront developments. For this reason, we have decided to continue developing openCFS (before 2020, known as CFS++ Coupled Field Simulations written in C++).

Finally, we select two acoustic simulations to show the acoustic capabilities of openCFS. The aeroacoustic simulation model simVoice for human phonation \cite{dollinger2023overview} was a 3D finite volume, finite element coupled model validated by a benchmarking dataset \cite{schoder2024benchmark}, and a nice summary is presented in \cite{falk20213d}. The unique feature of human phonation simulations is that the flow-induced sources must be captured precisely on a time-dependent geometry to get a good prediction quality.

\begin{figure}[ht!]
    \centering
    \includegraphics[width=0.5\linewidth]{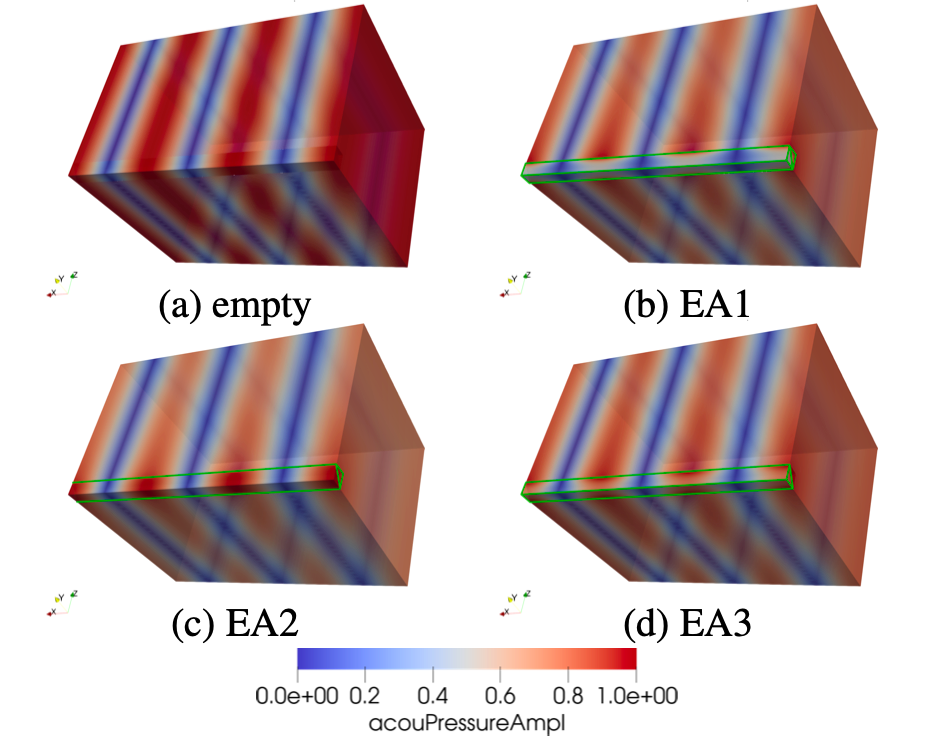}
    \caption{Results of the eigenvector for different configurations of the edge absorbers (highlighted in green) compared to the empty room configuration \cite{kraxberger2023nonlinear}.}
    \label{fig:nonlinear}
\end{figure}

Figure \ref{fig:nonlinear} shows the influence of edge absorbers on the nonlinear mode shapes of the acoustic field.
For this simulation, openCFS was used as the physical solver, and pyCFS triggered the iteration step for the fixed-point method.

\subsection{Outlook for 2025}
As an outlook for 2025, the frequency dependent material will be available in the master branch for time domain simulations using auxiliary variables. Furthermore, an automated generation of a PML domain and curvilinear PML is available. The capabilities of the mortar methods will be enhanced in acoustics and the acousticPDE class will be cleaned up.


%
\section{Acknowledge}
We would like to acknowledge the authors of openCFS:
Angermeier, Katharina; Bahr, Ludwig; Dev, Chaitanya; Eiser, Sebastian; Escobar, Max; Freidhager, Clemens; Grabinger, Jens; Greifenstein, Jannis; Guess, Thomas; Hassanpour Guilvaiee, Hamideh; Hauck, Andreas; Hofer, Fred; Hübner, Daniel; Hüppe, Andreas; Jaganathan, Srikrishna; Junger, Clemens; Kaltenbacher, Manferd; Landes, Herman; Link, Gerhard; Mayrhofer, Dominik; Michalke, Simon; Mohr, Markus; Nierla, Michael; Perchtold, Dominik; Roppert, Klaus; Schmidt, Bastian; Schoder, Stefan; Schury, Fabian; Seebacher, Philipp; Shaposhnikov, Kirill; Tautz, Matthias; Toth, Florian; Triebenbacher, Simon; Volk, Adian; Vu, Bich Ngoc; Wein, Fabian; Zhelezina, Elena; Zörner, Stefan.

Sketches in this work have been partly created using the Adobe Illustrator plug-in LaTeX2AI (\url{https://github.com/stoani89/LaTeX2AI}).

\bibliographystyle{plain}
\bibliography{references}  






\end{document}